\documentclass[12pt]{article}

\usepackage{amssymb}
\usepackage{amsmath}
\usepackage{xypic}

\begin{document}

\title{Probl\`eme de Lehmer pour les hypersurfaces de vari\'et\'es ab\'eliennes de type C.M.}
\author{Nicolas Ratazzi,\\
Universit\'e Paris 6, Projet th\'eorie des nombres, UMR 7586,\\
 case 247, 4 place Jussieu, Institut de math\'ematiques,\\
 75252 Paris, FRANCE}

\newcounter{ndefinition} 
\newcommand{\defi}{\addtocounter{ndefinition}{1}{\noindent \textbf{D{\'e}finition \thendefinition\ }}}
\newcounter{nrem}
\newcommand{\rem}{\addtocounter{nrem}{1}{\noindent \textbf{Remarque \thenrem\ }}}
\newtheorem{lemme}{Lemme} 
\newtheorem{coro}{Corollaire} 
\newtheorem{theo}{Th{\'e}or{\`e}me}
\newtheorem{conj}{Conjecture}
\newtheorem{prop}{Proposition}

\maketitle

\renewcommand{\thefootnote}{}
\footnote{\textit{Email address : }ratazzi@math.jussieu.fr}

\vspace{1cm}

\section{Introduction}
\noindent On sait depuis les travaux de Philippon \cite{phi1} \cite{phi2} \cite{phi3}, puis Bost, Gillet, Soul\'e \cite{BGS} dans le cadre de l'intersection arithm\'etique, comment d\'efinir la hauteur des vari\'et\'es projectives ; l'id\'ee \'etant de consid\'erer un point comme une vari\'et\'e de dimension z\'ero et de g\'en\'eraliser ceci en dimension sup\'erieure. De m\^eme que dans le cas des points, on sait pour les vari\'et\'es ab\'eliennes munies d'un fibr\'e en droites ample et sym\'etrique d\'efinir une hauteur particuli\`erement agr\'eable~: la hauteur canonique $\widehat{h}_L$, ou hauteur normalis\'ee. En dimension z\'ero, il existe un th\'eor\`eme caract\'erisant les points de hauteur normalis\'ee nulle ; c'est un r\'esultat de Kronecker dans le cas de $\mathbb{G}_m$. \mbox{Philippon} \cite{phi3} (dans le cas d'un produit de courbes elliptiques) puis Zhang \cite{zhang2} et David-Philippon \cite{daviphi} dans le cas g\'en\'eral ont montr\'e comment g\'en\'eraliser ce r\'esultat pour caract\'eriser les sous-vari\'et\'es de hauteur normalis\'ee nulle : ce sont les translat\'ees d'une sous-vari\'et\'e ab\'elienne par un point de torsion. On dit qu'une telle sous-vari\'et\'e est une sous-vari\'et\'e de torsion. La r\'eponse \`a cette question r\'esoud \`a une conjecture de Bogomolov. Ceci \'etant, on peut se demander comment minorer la hauteur normalis\'ee d'une sous-vari\'et\'e de hauteur non-nulle d'une vari\'et\'e ab\'elienne. Dans leur article \cite{daviphi}, David et Philippon ont formul\'e un probl\`eme g\'en\'eral (le probl\`eme 1.7) contenant cette question. En terme du degr\'e d\'efini ci-dessous, on peut notamment faire ressortir de la discussion suivant la formulation de leur probl\`eme l'\'enonc\'e suivant :

\begin{conj}\label{conj1} (David-Philippon) Soit $A$ une vari\'et\'e ab\'elienne d\'efinie sur un corps de nombres $k$, munie d'un fibr\'e ample et sym\'etrique $\mathcal{L}$. Soit $V$ une sous-vari\'et\'e stricte de $A$ sur $k$, $k$-irr\'eductible et qui n'est pas r\'eunion de sous-vari\'et\'es de torsion, alors, on a l'in\'egalit\'e
\[\frac{\widehat{h}_{\mathcal{L}}(V)}{\deg_{\mathcal{L}}(V)}\geq c(A,\mathcal{L})\deg_{\mathcal{L}}(V)^{-\frac{1}{s-\textnormal{dim} V}},\]
\noindent o\`u $s$ est la dimension du plus petit sous-groupe alg\'ebrique contenant $V$, et o\`u $c(A,\mathcal{L})$ est une constante ne d\'ependant que de $A$ et de $\mathcal{L}$.
\end{conj}

\subsection{Degr\'e et hauteur}
\noindent Soient $k$ un corps de nombres suppos\'e plong\'e dans $\mathbb{C}$, et $\mathcal{O}_k$ son anneau d'entiers. On dira que $V$ est une \textit{vari\'et\'e alg\'ebrique} sur $k$ si $V$ est un $k$-sch\'ema de type fini g\'eom\'etriquement r\'eduit. On dira que $G$ est un \textit{groupe alg\'ebrique} sur $k$ si c'est une vari\'et\'e en groupes sur $k$. On dira que $A$ est une \textit{vari\'et\'e ab\'elienne} d\'efinie sur $k$ si c'est un groupe alg\'ebrique connexe propre et lisse sur $k$. Par sous-vari\'et\'e on entendra toujours sous-vari\'et\'e ferm\'ee.

\vspace{.3cm}

\defi On dit qu'une vari\'et\'e ab\'elienne simple $A/k$ sur un corps de nombres est \textit{de type C.M.} si son anneau d'endomorphismes tensoris\'e par $\mathbb{Q}$ contient (apr\`es \'eventuellement extension du corps de base) un corps commutatif $F$ de dimension $2\,\textnormal{dim} A$ sur $k$. Une vari\'et\'e ab\'elienne $A/k$ est dite \textit{de type C.M.} si son anneau d'endomorphismes tensoris\'e par $\mathbb{Q}$ contient un produit de corps de nombres $K_1\times\cdots\times K_r$ tels que $\sum[K_i :\mathbb{Q}]=2\,\textnormal{dim} A$.

\vspace{.3cm}

\noindent Soit $X$ une vari\'et\'e projective munie d'un plongement $\varphi_{\mathcal{L}} : X \hookrightarrow \mathbb{P}^n_k$ d\'efini par un fibr\'e $\mathcal{L}$ tr\`es ample sur $X$. Si $\mathcal{O}(1)$ d\'enote le fibr\'e standard sur $\mathbb{P}^n_{\mathcal{O}_k},$ on a $\varphi_{\mathcal{L}}^*\mathcal{O}(1)_k\simeq \mathcal{L}$. On note $\overline{\mathcal{O}(1)}$ le fibr\'e standard muni de la m\'etrique de Fubini-Study. Si $V$ est une sous-vari\'et\'e de $X$, on note $\mathcal{V}_{\mathcal{L}}$ l'adh\'erence sch\'ematique de $\varphi_{\mathcal{L}}(V)$ dans $\mathbb{P}^n_{\mathcal{O}_k}$.

\vspace{.3cm}

\defi Si $\mathcal{L}$ est un fibr\'e ample sur une vari\'et\'e ab\'elienne $A$, et $V$ une sous-vari\'et\'e de $A$, on d\'efinit le \textit{degr\'e de la vari\'et\'e $V$}  relativement \`a $\mathcal{L}$, et on note $\deg_{\mathcal{L}} V$ l'entier 
$\deg \left(c_1(\mathcal{L})^{\textnormal{dim}V}\cdot V\right)$ o\`u $\deg$ est le degr\'e projectif usuel d'un $0$-cycle.

\vspace{.3cm}

\defi On appelle \textit{hauteur de la vari\'et\'e $V$} associ\'ee \`a $\mathcal{L}$, et on note $h_{\mathcal{L}}(V)$ le r\'eel $h_{\overline{\mathcal{O}(1)}}(\mathcal{V}_{\mathcal{L}})$ o\`u $h_{\overline{\mathcal{O}(1)}}(.)$ est la hauteur, au sens de Bost-Gillet-Soul\'e \cite{BGS}, associ\'ee au fibr\'e hermitien $\overline{\mathcal{O}(1)}$.

\vspace{.3cm}

\rem Par le th\'eor\`eme 3 p. 366 de \cite{soule}, $h_{\mathcal{L}}(V)$ coincide avec la hauteur $h(f_{V,\mathcal{L}})$ de Philippon, telle que d\'efinie au paragraphe 2. de \cite{phi3}, o\`u $f_{V,\mathcal{L}}$ est une forme \'eliminante de l'id\'eal de d\'efinition de $\varphi_{\mathcal{L}}(V)$ dans $k[X_0,\ldots,X_n]$. (Le terme d'erreur de \cite{soule} disparait du fait du changement de normalisation pour la hauteur de Philippon entre les articles \cite{phi1} et \cite{phi3}).

\vspace{.3cm}

\defi Dans le cas o\`u $X=A$ est une vari\'et\'e ab\'elienne, et o\`u $\mathcal{L}$ est de plus sym\'etrique, Philippon \cite{phi3}, puis Zhang \cite{zhang2} avec des m\'ethodes arakeloviennes, ont montr\'e en utilisant un proc\'ed\'e de limite \`a la N\'eron-Tate, comment d\'efinir une \textit{hauteur canonique}, not\'ee $\widehat{h}_{\mathcal{L}}(.)$, sur l'ensemble des sous-vari\'et\'es de $A$. Cette hauteur v\'erifie notamment : si $X$ est une sous-vari\'et\'e de $A$, de stabilisateur $G_X$, et si $n$ est un entier, alors,
\[\widehat{h}_{\mathcal{L}}\left([n](X)\right)=\frac{n^{2(\textnormal{dim} X+1)}}{\mid \textnormal{ker}\ [n]\cap G_X\mid} \widehat{h}_{\mathcal{L}}(X).\] 

\vspace{.3cm}

\defi Soit $A/k$ une vari\'et\'e ab\'elienne. On dit que $V$ est une \textit{sous-vari\'et\'e de torsion} de $A$ si $V=a+B$ avec $a\in A_{\textnormal{tors}}$ et $B$ une sous-vari\'et\'e ab\'elienne de $A$.

\vspace{.3cm}

\noindent D'apr\`es les r\'esultats de Philippon \cite{phi3}, David-Philippon \cite{daviphi} et Zhang \cite{zhang2}, on a, si $V$ est une sous-vari\'et\'e de $A/k$ d\'efinie sur une extension finie $K/k$,
\[\widehat{h}_{\mathcal{L}}(V)=0\textnormal{ si et seulement si $V$ est une sous-vari\'et\'e de torsion.}\]

\vspace{.3cm}

\defi Soient $V$ une sous-vari\'et\'e de $A$ sur $k$, et $\theta$ un nombre r\'eel positif. On pose $V(\theta,\mathcal{L})=\left\{x\in V(\overline{k})\ / \ \widehat{h}_{\mathcal{L}}(x)\leq \theta\right\}.$  On d\'efinit alors le \textit{minimum essentiel} de $V$, et on note $\hat{\mu}^{\textnormal{ess}}_{\mathcal{L}}(V)$ le r\'eel 
\[\hat{\mu}^{\textnormal{ess}}_{\mathcal{L}}(V)=\inf \left\{\theta>0 \ / \ \overline{V(\theta,\mathcal{L})}=V\ \right\}.\]

\subsection{R\'esultats}

\noindent Dans la direction de la conjecture \ref{conj1}, on a le r\'esultat suivant (cf. corollaire 2  de \cite{ratazzi})~:

\begin{theo}\label{theo1} Si $A$ est une vari\'et\'e ab\'elienne de type C.M., $\mathcal{L}$ un fibr\'e en droites ample et sym\'etrique de $A$, et si $V$ est une sous-vari\'et\'e alg\'ebrique stricte de $A$ sur $k$, $k$-irr\'eductible et qui n'est pas r\'eunion de sous-vari\'et\'es de torsion, alors, on a l'in\'egalit\'e
\[\frac{\widehat{h}_{\mathcal{L}}(V)}{\deg_{\mathcal{L}}(V)}\geq \hat{\mu}^{\textnormal{ess}}_{\mathcal{L}}(V)\geq c(A,\mathcal{L})\deg_{\mathcal{L}}(V)^{-\frac{1}{n-\textnormal{dim} V}}\left(\log(3\deg_{\mathcal{L}}(V))\right)^{-\kappa(n)},\]
\noindent o\`u $n$ est la dimension du plus petit sous-groupe alg\'ebrique contenant $V$, et o\`u $\kappa(n)$ est une constante effectivement calculable ne d\'ependant que de $n$ (par exemple $\kappa(n)=(2n(n+1)!)^{n+2}$ convient).
\end{theo}

\noindent On se restreint dans cet article au cas particulier des hypersurfaces $V$ d'une vari\'et\'e ab\'elienne de type C.M. Dans ce cas et sous les hypoth\`eses du th\'eor\`eme pr\'ec\'edent, on a n\'ecessairement $n=g$. En effet, par d\'efinition $n$ appartient \`a $\{g-1,g\}$. De plus, si $n$ \'etait \'egal \`a $g-1$, alors $V$ serait une r\'eunion de sous-vari\'et\'es de torsion, ce qui contredit l'hypoth\`ese faite sur $V$. Ainsi, dans le cas des hypersurfaces, la conjecture est la suivante :

\begin{conj}\label{conj2} Sous les hypoth\`eses pr\'ec\'edentes, et en supposant de plus que $V$ est hypersurface de $A$, on a l'in\'egalit\'e
\[\widehat{h}_{\mathcal{L}}(V)\geq c(A,\mathcal{L}),\]
\noindent o\`u $c(A,\mathcal{L})$ est une constante ne d\'ependant que de $A$ et de $\mathcal{L}$.
\end{conj}

\noindent De m\^eme, le th\'eor\`eme \ref{theo1} se sp\'ecialise en 

\begin{theo} Sous les hypoth\`eses pr\'ec\'edentes, et en supposant que $V$ est une hypersurface de $A$, on a l'in\'egalit\'e
\[\widehat{h}_{\mathcal{L}}(V)\geq \deg_{\mathcal{L}}(V)\hat{\mu}^{\textnormal{ess}}_{\mathcal{L}}(V)\geq c(A,\mathcal{L})\left(\log(3\deg_{\mathcal{L}}(V))\right)^{-\kappa(g)},\]
\noindent o\`u $g$ est la dimension de $A$, et o\`u $\kappa(g)$ est une constante effectivement calculable ne d\'ependant que de $g$ (par exemple $\kappa(g)=(2g(g+1)!)^{g+2}$ convient).
\end{theo}

\noindent Dans ce cadre restreint aux hypersurfaces, on montre un r\'esultat sensiblement plus fin en direction de la conjecture \ref{conj2}~: on peut prendre pour $\kappa$ une valeur absolue, ind\'ependante de $g$. En notant $\delta_{i,j}$ le symbole de Kronecker (valant $1$ si $i=j$ et $0$ sinon), on d\'emontre ici le r\'esultat suivant :

\begin{theo} Si $A$ est une vari\'et\'e ab\'elienne de type C.M., $\mathcal{L}$ un fibr\'e en droites ample et sym\'etrique de $A$, et si $V$ est une hypersurface irr\'eductible de $A$ sur $k$ qui n'est pas r\'eunion de sous-vari\'et\'es de torsion, alors, on a l'in\'egalit\'e
\[\widehat{h}_{\mathcal{L}}(V)\geq \deg_{\mathcal{L}}(V)\hat{\mu}^{\textnormal{ess}}_{\mathcal{L}}(V)\geq c(A,\mathcal{L})\frac{\left(\log\log \deg_{\mathcal{L}}V\right)^{1+2\delta_{g-s,1}}}{\left(\log\deg_{\mathcal{L}} V\right)^{2+\delta_{g-s,1}}},\]
\noindent o\`u $s$ est la dimension du stabilisateur de $V$.
\end{theo}

\noindent Notons que $\delta_{g-s,1}=0$ sauf si $A/k$ est le produit $E\times B$ d'une courbe elliptique $E/k$ et d'une vari\'et\'e ab\'elienne $B/k$, et si $V$ est de la forme $\overline{\{P\}}\times B$, o\`u $P$ est un point $\overline{k}$-rationnel de $E$ qui n'est pas de torsion. Dans ce cas, en supposant que $A/k$ est une courbe elliptique, $\mathcal{L}$ le fibr\'e associ\'e au diviseur $3(0)$, et o\`u $V=\overline{\{P\}}$ est l'ensemble des conjugu\'es d'un point non de torsion $P\in A(K)$ dans une extension finie $D=[K:k]$, on retrouve exactement le r\'esultat de Laurent \cite{laurent} sur le probl\`eme de Lehmer elliptique, \`a savoir
\[\widehat{h}(P)\geq \frac{c(A)}{D}\left(\frac{\log\log D}{\log D}\right)^3.\]
\noindent Dans le cas d'une ``vraie'' hypersurface, (i.e., quand $\delta_{g-s,1}=0$), on obtient une minoration un peu meilleure.

\vspace{.3cm}

\noindent La d\'emonstration suit fondamentalement les id\'ees (et reprend une grande partie des preuves) de l'article de David-Hindry \cite{davidhindry} concernant le probl\`eme de Lehmer pour les points d'une vari\'et\'e ab\'elienne. Il s'agit en fait d'une extension d'un travail de Amoroso-David \cite{amodavid2} concernant le cas des tores, au cas des vari\'et\'es ab\'eliennes de type C.M. On fait un raisonnement par l'absurde, et on se fixe une hypersurface $V$ contredisant la conclusion du th\'eor\`eme. La preuve consiste essentiellement en une preuve de transcendance classique. On commence tout d'abord par construire une fonction auxiliaire, nulle avec un grand ordre sur $V$. Pour cela on met en oeuvre une astuce d\^ue \`a Amoroso-David (qu'ils introduisent dans \cite{amodavid2}) permettant de se ramener \`a un syst\`eme d'\'equations fini et de hauteur control\'ee. Ceci nous permet d'appliquer un lemme de Siegel pour construire la fonction auxiliaire $F$. La deuxi\`eme partie de la preuve, l'extrapolation, consiste \`a montrer que $F$ continue \`a s'annuler avec un ordre relativement grand sur les transform\'ees $\alpha_v(V)$ de $V$ par certaines isog\'enies $\alpha_v$ o\`u $v$ d\'ecrit un ensemble de places finies convenables du corps de d\'efinition $k$ (les $\alpha_v$ sont des relev\'ees sur $A/k$ des morphismes de Frobenius en caract\'eristique finie $p_v$. C'est pour assurer l'existence de ces isog\'enies que l'on se restreint au cas C.M.). Il s'agit d'une extrapolation aux places $v$-adiques. L'id\'ee pour montrer ceci est d'appliquer une g\'en\'eralisation du petit th\'eor\`eme de Fermat : c'est la m\'ethode employ\'ee pour la premi\`ere fois par Dobrowolski \cite{dob} dans le cas du probl\`eme de Lehmer sur $\mathbb{G}_m$. Cette id\'ee a ensuite \'et\'e reprise par Laurent \cite{laurent} dans le cas des courbes elliptiques \`a multiplication complexes puis \'etendue au cas des vari\'et\'es ab\'eliennes de type C.M. par David-Hindry \cite{davidhindry}. C'est cette derni\`ere g\'en\'eralisation que nous allons reprendre. Ceci \'etant fait, il suffit pour conclure d'appliquer le th\'eor\`eme de B\'ezout g\'eom\'etrique pour aboutir \`a une contradiction (pour peu que les diff\'erents param\`etres intervenant dans l'\'etape de transcendance aient \'et\'es convenablement choisis). Pour cette derni\`ere \'etape, on a besoin d'avoir une bonne minoration du degr\'e de l'union des $\alpha_v(V)$. Ceci se fait en suivant les calculs de $\cite{davidhindry}$.

\vspace{.3cm}

\noindent \textbf{Remerciements :} Je tiens \`a remercier Sinnou David pour m'avoir sugg\'erer l'\'ecriture de cet article, et je tiens \'egalement \`a remercier Marc Hindry pour les nombreuses discussions que nous avons eu sur le sujet.

\section{Frobenius, isog\'enies admissibles et d\'e\-ri\-va\-tions}
\subsection{Morphismes de Frobenius}
\noindent On commence par introduire quelques notations : 

\vspace{.3cm}

\noindent Si $k$ est un corps de nombres, on note $\mathcal{O}_k$ son anneau d'entiers, $v$ une place finie de $k$, et $k_v$ le corps r\'esiduel associ\'e \`a $v$.

\vspace{.3cm}

\noindent Si $A/k$ est une vari\'et\'e ab\'elienne, on note $\mathcal{A}/\mathcal{O}_k$ son mod\`ele de N\'eron, et $A_v/k_v$ la fibre sp\'eciale correspondant \`a la place finie $v$. Rappelons la propri\'et\'e universelle du mod\`ele de N\'eron : si $\mathcal{X}/\mathcal{O}_k$ est lisse, de fibre g\'en\'erique $X/k$, tout $k$-morphisme $X\rightarrow A$ se rel\`eve de mani\`ere unique en un $\mathcal{O}_k$-morphisme $\mathcal{X}\rightarrow \mathcal{A}$.

\vspace{.3cm}

\noindent Sur la vari\'et\'e $A_v/k_v$, on dispose d'un endomorphisme particulier : le morphisme de Frobenius $\textnormal{Frob}_v$, correspondant en coordonn\'ees projectives \`a l'\'el\'e\-va\-tion \`a la puissance $q=\textnormal{N}(v)$, o\`u $\textnormal{N}(v)$ est la norme $K/\mathbb{Q}$ de $v$.

\vspace{.3cm}

\noindent La propri\'et\'e universelle du produit fibr\'e $A_v=\mathcal{A}\times_{\mathcal{O}_k}k_v$ permet d'associer naturellement \`a tout $\mathcal{O}_k$-endomorphisme de $\mathcal{A}$ un $k_v$-endomorphisme de $A_v$. En utilisant la propri\'et\'e universelle du mod\`ele de N\'eron, on en d\'eduit une fl\`eche naturelle 
\[\Psi : \textnormal{End}_k(A) \rightarrow \textnormal{End}_{k_v}(A_v).\]
Cette fl\`eche n'est en g\'en\'eral pas surjective, mais on peut par contre montrer qu'elle est injective aux places de bonne r\'eduction. Dans le cas C.M., un th\'eor\`eme de Shimura-Taniyama permet d'affirmer que le morphisme $\textnormal{Frob}_v$ se rel\`eve en presque toutes places :

\begin{prop}\label{shimura}\textnormal{(Shimura-Taniyama)} Soit $A/k$ une vari\'et\'e ab\'elienne de ty\-pe C.M. Notons $\prod_{i=1}^r K_i$ le produit de corps de nombres qui est inclus dans $\textnormal{End}_k(A)\otimes\mathbb{Q}$ et tel que $\sum_{i=1}^r[K_i :\mathbb{Q}]=2\,\textnormal{dim} A$. On suppose que le corps de nombres $k$ contient tous les $K_i$, et que $\prod_{i=1}^r\mathcal{O}_{K_i}$ est inclus dans $\textnormal{End}_k(A)$. Alors, pour presque toutes places, l'endomorphisme $\textnormal{Frob}_v$ se rel\`eve en un $k$-endomorphisme $\alpha_v$ de $A$. On appelera morphisme de Frobenius sur $A$ un tel endomorphisme.
\end{prop}
\noindent \textbf{D\'emonstration} C'est le Theorem 1 paragraphe III.13 de \cite{shimura}.\hfill$\Box$

\vspace{.3cm}

\noindent Ce sont ces morphismes de Frobenius sur $A/k$ qui vont nous permettre d'\'ecrire l'\'etape d'extrapolation.

\vspace{.3cm}

\rem En fait on pourrait sp\'ecifier les places qu'il faut exclure dans la proposition, mais nous n'en aurons pas besoin. Par ailleurs, pour pouvoir appliquer le th\'eor\`eme, il faut v\'erifier deux conditions : la premi\`ere est toujours satisfaite quitte \`a faire une extension de degr\'e born\'e de $k$. La seconde n'est pas toujours satisfaite, mais on peut toujours trouver une vari\'et\'e ab\'elienne isog\`ene qui la v\'erifie. 

\vspace{.3cm}

\noindent Quitte \`a faire une extension de degr\'e born\'e de $k$, et quitte \`a prendre une vari\'et\'e ab\'elienne isog\`ene \`a la vari\'et\'e de d\'epart, on supposera d\'esormais toujours que les hypoth\`eses de la proposition \ref{shimura} sont satisfaites.

\subsection{Isog\'enies admissibles}
\noindent On rappelle la notion d'isog\'enie admissible telle qu'introduite dans \cite{davidhindry}.

\vspace{.3cm}

\defi Soient $A$ une vari\'et\'e ab\'elienne et $\mathcal{L}$ un fibr\'e ample sur $A$. Une isog\'enie $\alpha$ de $A$ est dite \textit{admissible} par rapport \`a $\mathcal{L}$ si 
\begin{enumerate}
\item $\alpha$ est dans le centre de $\textnormal{End}(A)$.

\vspace{.3cm}

\item il existe un entier $\textnormal{q}(\alpha)$ appel\'e \textit{poids} de $\alpha$ tel que $\alpha^{\star}\mathcal{L}\simeq \mathcal{L}^{\otimes \, \textnormal{q}(\alpha)}$.
\end{enumerate}

\vspace{.3cm}

\rem En fait la condition (1) ne sert qu'a simplifier l'\'enonc\'e du lemme \ref{distinct}. C'est la condition (2) qui importe vraiment. Les seules isog\'enies qui nous int\'eresseront sont les relev\'ees $\alpha_v$ des morphismes de Frobenius qui sont admissibles (cf. la Proposition \ref{frob}).

\begin{lemme}\label{degre} Soient $A$ une vari\'et\'e ab\'elienne de dimension $g$ munie d'un fibr\'e en droites tr\`es ample $\mathcal{L}$, et $\alpha$ une isog\'enie admissible relativement \`a $\mathcal{L}$, de poids $q=\textnormal{q}(\alpha)$. Dans le plongement projectif de $A$, associ\'e \`a $\mathcal{L}$, $A\hookrightarrow \mathbb{P}_n$, on a : 
\begin{enumerate}
\item $\textnormal{card}\left(\textnormal{ker}(\alpha)\right)=q^g$,
\item pour toute sous-vari\'et\'e $V$ de $A$ de stabilisateur $G_V$, on a 
\[\deg_{\mathcal{L}}\left(\alpha(V)\right)=\frac{q^{\dim(V)}}{\left| G_V\cap \textnormal{ker}(\alpha)\right|}\deg_{\mathcal{L}}(V)\]
\end{enumerate}
\end{lemme}
\noindent \textbf{D\'emonstration} Le point (1) est facile : par d\'efinition, $\alpha^{\star}\mathcal{L}\simeq \mathcal{L}^{\otimes q}$. On a donc, 
\[q^g\deg_{\mathcal{L}}(A)=\deg_{\mathcal{L}^{\otimes q}}(A)=\deg_{\alpha^{\star}\mathcal{L}}(A)=\left|\textnormal{ker}(\alpha)\right|\deg_{\mathcal{L}}(A).\]
\noindent L'amplitude de $\mathcal{L}$ nous assure que le dernier degr\'e est strictement positif. On simplifie pour conclure. Pour le point (2), il s'agit du point (ii) du lemme 6. de \cite{hindry}.\hfill$\Box$

\vspace{.3cm}

\begin{lemme}\label{lemme25}Soient $G$ un sous-groupe alg\'ebrique de la vari\'et\'e ab\'elienne $A/k$, $\mathcal{L}$ un fibr\'e tr\`es ample sur $A$, et $\alpha$ une isog\'enie admissible relativement \`a $\mathcal{L}$ de poids $\textnormal{q}(\alpha)$ de $A$. On a 
\[ \textnormal{q}(\alpha)^{\dim G}\leq\textnormal{card}\left(\textnormal{ker}(\alpha)\cap G\right)\leq \left[G:G^0\right]\textnormal{q}(\alpha)^{\dim G}.\]
\end{lemme}
\noindent \textbf{D\'emonstration} On note que 
\[ \left[G:G^0\right]\textnormal{card}\left(\textnormal{ker}(\alpha)\cap G^0\right)\geq \textnormal{card}\left(\textnormal{ker}(\alpha)\cap G\right)\geq \textnormal{card}\left(\textnormal{ker}(\alpha)\cap G^0\right).\]
\noindent La restriction de $\alpha$ \`a la sous-vari\'et\'e ab\'elienne $G^0$ est encore une isog\'enie admissible de poids $\textnormal{q}(\alpha)$ pour $(G^0,\mathcal{L}_{\mid G^0})$ (cf. Lemme 2.4. point (ii) de \cite{davidhindry}). Par le point (1) du lemme \ref{degre} pr\'ec\'edent, on en d\'eduit que le cardinal du noyau de cette isog\'enie $\alpha_{\mid G^0}$ est $\textnormal{q}(\alpha)^{\dim G^0}$. \hfill $\Box$

\vspace{.3cm}

\noindent Soit $V$ une sous-$k$-vari\'et\'e stricte de $A$, $k$-irr\'eductible. Le lemme suivant (dont l'origine remonte \`a Dobrowolski \cite{dob}) montre que les images par une isog\'enie admissible de ses composantes g\'eom\'etriquement irr\'eductibles sont essentiellement distinctes. On commence pour cela par donner une d\'efinition~:

\vspace{.3cm}

\defi Soient $A$ une vari\'et\'e ab\'elienne et $\mathcal{L}$ un fibr\'e en droites ample sur $A$. Deux isog\'enies admissibles de $A$ par rapport \`a $\mathcal{L}$ sont dites \textit{premi\`eres entre elles} si leurs poids sont premiers entre eux.

\begin{lemme}\label{distinct} Soient $A$ une vari\'et\'e ab\'elienne sur $k$ de dimension $g\geq 1$, $\mathcal{L}$ un fibr\'e en droites tr\`es ample sur $A$, $V$ une sous-$k$-vari\'et\'e stricte de $A$, irr\'eductible sur $k$.  Si $V$ n'est pas une r\'eunion de sous-vari\'et\'es de torsion de $A$, on a :
\begin{enumerate}
\item Pour tout couple $(\alpha,\beta)$ d'isog\'enies admissibles pour $\mathcal{L}$, de poids distincts, pour tout $\sigma\in \textnormal{Gal}(\overline{k}/k)$, et pour toute composante g\'eom\'etri\-quement irr\'e\-ductible $W$ de $V$, les sous-vari\'et\'es $\alpha(W)$ et $\beta\left(\sigma(W)\right)$ sont distinctes.

\vspace{.3cm}

\item Soit $\mathcal{P}$ un ensemble d'isog\'enies admissibles pour $\mathcal{L}$, deux \`a deux pre\-mi\`eres entre elles. Notons $V_1,\ldots, V_M$ les composantes g\'eom\'etri\-que\-ment irr\'e\-duc\-tibles de $V$, et notons $\mathcal{Q}$ le sous-ensemble de $\mathcal{P}$ d\'efini par
\[ \mathcal{Q}=\left\{ \alpha\in \mathcal{P}\ / \ \exists i,j,\ 1\leq i<j\leq M,\ \ \alpha(V_i)=\alpha(V_j)\right\}.\]
\noindent Le cardinal de $\mathcal{Q}$ est major\'e par $\frac{\log M}{\log 2}$.
\end{enumerate}
\end{lemme}
\noindent \textbf{D\'emonstration} Dans ce contexte il s'agit de la proposition 2.7. de \cite{davidhindry}\hfill$\Box$

\vspace{.3cm}

\noindent On conclut ce paragraphe en ``rappelant'' que les morphismes de Frobenius sur $A/k$ sont des isog\'enies admissibles :

\vspace{.3cm}

\defi Soient $A$ une vari\'et\'e ab\'elienne et $\mathcal{L}$ un fibr\'e en droites ample sur $A$. Suivant Mumford , on dit que $\mathcal{L}$ est \textit{totalement sym\'etrique} si $\mathcal{L}$ est le carr\'e d'un fibr\'e sym\'etrique.

\vspace{.3cm}

\noindent Le th\'eor\`eme de Lefschetz (cf. par exemple le Theorem A.5.3.6 de \cite{hindrysil}) nous indique que si $\mathcal{L}$ est un fibr\'e ample, alors $\mathcal{L}^{\otimes 3}$ est tr\`es ample.

\begin{prop}\label{frob}Soient $A/k$ une vari\'et\'e ab\'elienne de type C.M. v\'erifiant les hypoth\`eses de la proposition \ref{shimura}, et $\mathcal{L}$ un fibr\'e tr\`es ample et totalement sy\-m\'etri\-que sur $A$. Soit $\alpha_v$ un morphisme de Frobenius sur $A$ pour la place finie $v$. Alors, $\alpha_v$ est une isog\'enie admissible pour $\mathcal{L}$ de poids $\textnormal{q}(\alpha)$.
\end{prop}
\noindent \textbf{D\'emonstration} C'est la proposition 3.3. de \cite{davidhindry}.\hfill$\Box$

\section{Donn\'ees}
\subsection{Situation} 
\defi On dit qu'une sous-vari\'et\'e $X$ de $\mathbb{P}_n$ est \textit{projectivement normale} si son anneau de coordonn\'ees $S(X)$ est un anneau normal (i.e., int\'egralement clos).

\vspace{.3cm}

\noindent On peut montrer (cf. par exemple Birkenhake-Lange \cite{lange} p. 190-193) que $X\subset \mathbb{P}_n$ est projectivement normale si et seulement si elle est normale, et pour tout $d\geq 0$ la fl\`eche naturelle 
\[ H^0(\mathbb{P}_n,\mathcal{O}_{\mathbb{P}_n}(d))\rightarrow H^0(X,\mathcal{O}_X(d))\]
\noindent est surjective.

\vspace{.3cm}

\noindent Concernant les vari\'et\'es ab\'eliennes plong\'ees de mani\`ere projectivement normale, on a le r\'esultat suivant que l'on trouve dans \cite{lange} theorem 3.1 p. 190.

\begin{prop}Soient $A/k$ une vari\'et\'e ab\'elienne, et $\mathcal{L}$ un fibr\'e ample sur $A$. Pour tout $n\geq 3$, le fibr\'e $\mathcal{L}^{\otimes n}$ d\'efinit un plongement projectivement normal de $A$ dans un espace projectif $\mathbb{P}_n$.
\end{prop}

\vspace{.3cm}

\noindent Soient $A/k$ une vari\'et\'e ab\'elienne sur un corps de nombres, munie d'un fibr\'e sym\'etrique ample $\mathcal{L}$. Quitte \`a travailler avec $\mathcal{L}^{\otimes 4}$ plut\^ot qu'avec $\mathcal{L}$, on peut supposer que $\mathcal{L}$ est tr\`es ample, totalement sym\'etrique et d\'efinit un plongement projectivement normal de $A$ dans un projectif $\mathbb{P}_n$. On note $\mathcal{M}=\mathcal{L}\boxtimes\mathcal{L}$ le fibr\'e sur $A\times A$ associ\'e \`a $\mathcal{L}$. Soit $V$ une $k$-hypersurface irr\'eductible de $A$. On note $I_V$ l'id\'eal de d\'efinition de $V$ dans $\mathbb{P}_n$. Si $N$ est un entier, on a  
\[
\begin{array}{ccccccccc}
V & \subset 	& A & \overset{i}{\hookrightarrow} & A\times A  & \hookrightarrow & \mathbb{P}_n\times\mathbb{P}_n & \underset{\textnormal{Segre}}{\hookrightarrow} & \mathbb{P}_{(n+1)^2-1}\\
  & 		& x & \mapsto         & (x,[N]x) &                 &          &                                          &
\end{array}
\]

\vspace{.3cm}

\noindent Soient $L$ et $T$ deux entiers. On note $\left\{s_0,\ldots,s_l\right\}$ une base de $H^0(A\times A,\mathcal{M})$. On peut, par projective normalit\'e, choisir une base $\left\{Q_1,\ldots,Q_m\right\}$ du $k$-vectoriel $H^0\left(A\times A,\mathcal{M}^{\otimes L}\right)$ telle que tous les $Q_i$ sont homog\`enes de degr\'e $L$ en les $s_j$. De plus, on peut aussi voir les $s_i$ comme des $(1,1)$-formes homog\`enes de $k[\textbf{X},\textbf{Y}]$ o\`u $\textbf{X}=(X_0,\ldots,X_n)$, et $\textbf{Y}=(Y_0,\ldots,Y_n)$. Enfin on note $T_B$ l'espace tangent \`a l'origine de la sous-vari\'et\'e ab\'elienne $B=i(A)$ de $A\times A$ d\'efinie par $y=[N]x$. 

\subsection{Choix des param\`etres}

\noindent Soit $C_0$ un r\'eel positif, on note $s$ la dimension du stabilisateur de $V$, et $\delta_{i,j}$ le symbole de Kronecker (valant $1$ si $i=j$ et $0$ sinon). On  pose
\[ N_1=\left[C_0^{g+2}\left(\log\deg_{\mathcal{L}}V\right)^{1+\delta_{g-s,1}}\left(\log\log \deg_{\mathcal{L}}V\right)^{1-2\delta_{g-s,1}}\right],\]
\[ m=\left[\frac{\log\left(C_0^{\frac{g+1}{2}}\left(\deg_{\mathcal{L}}V\right)^{\frac{1}{2}}\left(\log\deg_{\mathcal{L}}V\right)^{\frac{1}{2}}\left(\log\log\deg_{\mathcal{L}}V\right)^{-1}\right)}{\log 2}\right],\  N=2^{m+1}\]
\[ T=\left[C_0^{g+1}\deg_{\mathcal{L}}V\log \deg_{\mathcal{L}}V\left(\log\log\deg_{\mathcal{L}}V\right)^{-3}\right],\]
\[  L=\left[C_0^{g+\frac{1}{2}}\deg_{\mathcal{L}}V\log \deg_{\mathcal{L}}V\left(\log\log\deg_{\mathcal{L}}V\right)^{-2}\right],\]
\noindent et,
\[ T_1=\left[C_0^{g}\deg_{\mathcal{L}} V\left(\log\log\deg_{\mathcal{L}}V\right)^{-2}\right].\]

\vspace{.3cm}

\noindent Ces param\`etres sont choisis de sorte que :
\begin{enumerate}
\item le nombre $N$ est une puissance de $2$ et v\'erifie l'encadrement
\[\frac{N}{2}\leq C_0^{\frac{g+1}{2}}\left(\deg_{\mathcal{L}}V\right)^{\frac{1}{2}}\left(\log\deg_{\mathcal{L}}V\right)^{\frac{1}{2}}\left(\log\log\deg_{\mathcal{L}}V\right)^{-1}< N.\]
\item $N^2>L+1$, afin qu'une forme $F$ bihomog\`ene de bi-degr\'e $(L,L)$ qui est non-identiquement nulle sur $A\times A$, ne soit pas identiquement nulle sur la sous-vari\'et\'e ab\'elienne $B$.
\item le minimum essentiel des vari\'et\'es intervenant est born\'e, autrement dit, 
\[N^2N_1\hat{\mu}^{\textnormal{ess}}(V)\leq c,\]
\item $T>L$, o\`u $T$ va \^etre l'ordre d'annulation dans le lemme de Siegel, et $L$ le degr\'e du polyn\^ome construit.
\item $T>T_1$, puisqu'on ne peut pas, par extrapolation esp\'erer un ordre d'annu\-la\-tion meilleur que celui dont on est parti ($T_1$ \'etant l'ordre d'annulation sur les sous-vari\'et\'es sur lesquelles on extrapole).
\end{enumerate}

\vspace{.3cm}

\noindent On fixe un premier $p_0$ (ne d\'ependant que de $A$) tel que pour tout premier $p\geq p_0$ et pour toute place $v$ divisant $p$, le morphisme de Frobenius $\alpha_v$ sur $A$ existe. On fixe alors pour chaque premier $p\geq p_0$ une place $v$ au dessus de $p$. On note $\mathcal{P}_k$ l'ensemble des places ainsi obtenues.

\vspace{.3cm}

\noindent Dans toute la suite, les in\'egalit\'es que l'on \'ecrira seront vraies pour tout $\deg_{\mathcal{L}}V$  et $C_0$ assez grands (i.e., plus grands qu'une constante ne d\'ependant que du couple $(A,\mathcal{L})$).

\section{Lemme de Siegel}

\noindent \textbf{But :} fabriquer un polyn\^ome, $F=\sum_{i=1}^m b_iQ_i$, \`a coefficients entiers relatifs, en les fonctions ab\'eliennes de $A\times A$, tel que $F$ est de ``petite'' hauteur, et tel que $F$ s'annule \`a un ordre sup\'erieur \`a $T$ sur $i(V)$, le long de $T_B$.

\vspace{.3cm}

\noindent En notant $\Theta$ l'application th\^eta d\'efinie sur $T_{A(\mathbb{C})}$ par la composition 
\[
\xymatrix{
T_{A(\mathbb{C})}\ \ar[r]^{\textnormal{exp}_{A(\mathbb{C})}}	& \  A(\mathbb{C}) \ar[r]^{\varphi_{\mathcal{L}}}	& \ \mathbb{P}_n }
\]
\noindent associ\'ee \`a $\mathcal{L}$, ceci correspond \`a trouver une solution de petite hauteur au syst\`eme d'inconnues les $b_i$
\begin{equation}\label{systeme1}
 \frac{\partial^{\kappa}F\left(\Theta(\textbf{u}+\textbf{z}),\Theta(N(\textbf{u}+\textbf{z})\right)}{\partial \textbf{z}^{\kappa}}_{\textnormal{\Large{$\mid$}}\textbf{z}=0}=0,
\end{equation}
\noindent pour tout $\mid \kappa\mid\leq T$ et $\textbf{u}\in T_{A(\mathbb{C})}$ tels que $\Theta(\textbf{u})\in V(\overline{k})$. 

\begin{lemme}\label{reduction1} Soit $\theta>\hat{\mu}^{\textnormal{ess}}_{\mathcal{L}}(V).$ Il existe un entier $d_0$ tel que si $F$ est une solution du syst\`eme
\begin{equation}\label{systeme2}
 \frac{\partial^{\kappa}F\left(\Theta(\textbf{u}+\textbf{z}),\Theta(N(\textbf{u}+\textbf{z})\right)}{\partial \textbf{z}^{\kappa}}_{\textnormal{\Large{$\mid$}}\textbf{z}=0}=0,
\end{equation}
\noindent pour tout $\mid \kappa\mid\leq T$ et $\textbf{u}\in T_{A(\mathbb{C})}$ tels que $\Theta(\textbf{u})$ appartient \`a l'ensemble fini 
\[ S_{d_0}(\theta)=\left\{x\in V(\overline{k})\ /\ \widehat{h}_{\mathcal{L}}(x)\leq \theta, \ \ [k(x): k]\leq d_0\right\},\]\noindent alors, $F$ est une solution du syst\`eme (\ref{systeme1}).
\end{lemme}
\noindent \textbf{D\'emonstration} Soit $d\geq 0$ un entier. On peut noter que l'ensemble $S_d(\theta)$ est stable sous l'action de $\textnormal{Gal}(\overline{k}/k)$ car $V$ est une $k$-vari\'et\'e. Par ailleurs, on a clairement $S_d(\theta)\subset S_{d+1}(\theta)$ pour tout $d\geq 0$. Notons $k[\textbf{X}]_L$ le $k$-espace vectoriel des polyn\^omes homog\`enes de degr\'e $L$, et $\mathcal{A}_d(\theta)$ le sous-$k$-espace vectoriel associ\'e \`a $S_d(\theta)$. La suite $\left(\mathcal{A}_d(\theta)\right)_{d\in \mathbb{N}}$ est une suite d\'ecroissante d'espaces vectoriels de dimension finie, elle est donc stationnaire. Notons $d_0$ l'indice \`a partir duquel cette suite est stationnaire. Par ailleurs, tous ces espaces contiennent le $k$-vectoriel $I_V^{(T)}{}_{\mid L}$ o\`u $I^{(T)}$ est la puissance symbolique $T$-i\`eme de $I$. Par d\'efinition de $d_0$, si $P$ est un polyn\^ome homog\`ene de degr\'e $L$ nul sur $S_{d_0}$, il appartient \`a $\mathcal{A}_{d_0}(\theta)$, et donc il s'annule sur $\bigcup_{d\geq 0} S_d(\theta)$. De plus, par d\'efinition du minimum essentiel, $\bigcup_{d\geq 0} S_d(\theta)$ est Zariski-dense dans $V$, donc le polyn\^ome $P$ s'annule sur $V$. Ainsi, dans le syst\`eme (\ref{systeme1}), on peut se restreindre aux $\textbf{u}\in T_{A(\mathbb{C})}$ tels que $\Theta(\textbf{u})$ appartient \`a $S_{d_0}(\theta)$. Par un th\'eor\`eme classique de Northcott, cet ensemble est fini.\hfill$\Box$

\vspace{.3cm}

\noindent On appelle syst\`eme (\ref{systeme2}) le nouveau syst\`eme ainsi obtenu. On passe maintenant \`a une estimation du rang.

\begin{lemme}\label{rang} Il existe une constante $c_1$ telle que le rang du syst\`eme (\ref{systeme1}) est major\'e par 
\[ c_1T(LN^2)^{g-1}\deg_{\mathcal{L}} V.                \]
\end{lemme}
\noindent \textbf{D\'emonstration} Il s'agit du lemme (ou plut\^ot de la preuve du lemme) 5.1 de \cite{davidhindry}. En effet, dans ce lemme, les auteurs de \cite{davidhindry} cherchent \`a obtenir une majoration du rang du syt\`eme
\begin{equation}
 \frac{\partial^{\kappa}F\left(\Theta(\textbf{u}+\textbf{z}),\Theta(N(\textbf{u}+\textbf{z})\right)}{\partial \textbf{z}^{\kappa}}_{\textnormal{\Large{$\mid$}}\textbf{z}=0}=0,
\end{equation}
o\`u $\textbf{u}$ est le logarithme d'un point $Q$ fix\'e. L'id\'ee est d'appliquer ``l'astuce de Philippon-Waldschmidt'' (voir \cite{phiwal} paragraphe 6, lemme 6.7).  Pour majorer le rang de ce syst\`eme, ils se donnent une vari\'et\'e $V$ de dimension $d$ contenant le point $Q$, et il majorent le syst\`eme
\begin{equation}
 \frac{\partial^{\kappa}F\left(\Theta(\textbf{u}+\textbf{z}),\Theta(N(\textbf{u}+\textbf{z})\right)}{\partial \textbf{z}^{\kappa}}_{\textnormal{\Large{$\mid$}}\textbf{z}=0}=0,
\end{equation}
\noindent pour tout $\mid \kappa\mid\leq T$ et $\textbf{u}\in T_{A(\mathbb{C})}$ tels que $\Theta(\textbf{u})\in V(\overline{k})$. Il obtiennent comme majorant du rang de ce syst\`eme le nombre  $c_1T^{g-d}(LN^2)^{d}\deg_{\mathcal{L}} V$. (on remplace dans leurs notations $T_0$ par $T$). En appliquant ceci \`a l'hypersurface $V$ consid\'er\'ee, on obtient donc le r\'esultat cherch\'e.\hfill$\Box$

\vspace{.3cm}

\noindent On peut maintenant \'enoncer le lemme de Siegel qui nous int\'eresse. Si $F=\sum a_{\textbf{i}}\textbf{X}^{\textbf{i}}$ est un polyn\^ome à coefficients dans $\overline{k}$, on d\'efinit classiquement sa hauteur $h(F)$ comme \'etant la hauteur logarithmique absolue du point projectif d\'efini par $1$ et tous les coefficients $a_i$ de $F$.

\vspace{.3cm}

\noindent L'objectif de l'article consiste \`a montrer que, $\deg_{\mathcal{L}}(V)\hat{\mu}^{\textnormal{ess}}_{\mathcal{L}}(V)>\frac{c(A,\mathcal{L})}{\log \deg_{\mathcal{L}}(V)^{\alpha}}$. On peut donc toujours supposer que $\hat{\mu}^{\textnormal{ess}}_{\mathcal{L}}(V)$ est strictement inf\'erieur \`a $1$.

\begin{prop}\label{siegel} Il existe une solution $F=\sum_{i=1}^m b_iQ_i$, $b_i\in \mathbb{Z}$ du syst\`eme (\ref{systeme1}) de degr\'e $L$ et de hauteur 
\[ h(F)\leq c_2C_0^{\frac{1}{2}(g+1)}\deg_{\mathcal{L}}(V)\log\deg_{\mathcal{L}}(V)(\log \log\deg_{\mathcal{L}} V)^{-2}.		\]
\end{prop}
\noindent \textbf{D\'emonstration} Soit $1>\theta>\hat{\mu}^{\textnormal{ess}}_{\mathcal{L}}(V).$ Par le lemme \ref{reduction1} il suffit, pour trouver une solution du syst\`eme (\ref{systeme1}), de trouver une solution du syst\`eme (2). Ceci remarqu\'e, on est ramen\'e \`a une preuve classique. On suit pour cela la preuve du lemme 5.4. de \cite{davidhindry}. 

\vspace{.3cm}

\noindent On commence par \'evaluer la hauteur de syst\`eme (2). Le syst\`eme (18) ainsi que l'in\'egalit\'e qui suit p. 42 de \cite{davidhindry} nous indique que la hauteur de chaque coefficient du syst\`eme est major\'ee par
\begin{equation}\label{systeme18}
c'_4LN^2\theta+T\left(\log(T+L)+\log N\right).
\end{equation}
\noindent Par ailleurs, le nombre d'inconnues $I$ est $\dim H^0\left(A\times A,\mathcal{M}^{\otimes L}\right)$. Le th\'eo\-r\`e\-me de Riemann-Roch pour les vari\'et\'es ab\'eliennes nous assure que $I$ v\'erifie l'en\-cadre\-ment
\begin{equation}\label{inconnues}
c_5L^{2g}\leq I\leq c_6L^{2g}.
\end{equation}
\noindent Notons $M$ la matrice du syst\`eme (2). Elle est d\'efinie sur $k$, donc si $\mathfrak{B}$ d\'enote le noyau de $M$, il est muni d'une $k$-structure. Si $\mathfrak{b}$ est la dimension de $\mathfrak{B}$, on a $\mathfrak{b}=I-\textnormal{rg}(M)$. Le lemme de Siegel classique (cf. par exemple Schmidt \cite{schmidt} Lemma IVB, p.10) nous indique alors qu'il existe une solution non-triviale comme recherch\'ee, de hauteur
\begin{equation}\label{hauteur}
h(F)\leq c_7\frac{h\left(\mathfrak{B}\right)}{\mathfrak{b}},
\end{equation}
o\`u $h\left(\mathfrak{B}\right)$ repr\'esente la hauteur du point $\mathfrak{B}$ d\'efini dans la grassmannienne correspondante. De plus, le Lemma IV p.10 de \cite{schmidt} nous indique que $h\left(\mathfrak{B}\right)=h\left(\mathfrak{B}^{\perp}\right).$ L'espace $\mathfrak{B}^{\perp}$ \'etant l'espace vectoriel engendr\'e par les colonnes de $M$, sa hauteur est par d\'efinition major\'ee par celle d'un mineur maximal $\Delta_{\max}$ de $M$. Cette derni\`ere hauteur est major\'ee par 
\begin{align*} 
h(\Delta_{\max})	& \leq  c_8\textnormal{rg}(M)\left(\log(\textnormal{rg} M)+c'_4LN^2\theta+T\left(\log(T+L)+\log N\right)\right)\\
			& \leq  c_9T(LN^2)^{g-1}\deg_{\mathcal{L}} V\left(\log(\deg_{\mathcal{L}} V)+2T\log T\right),
\end{align*}
\noindent la premi\`ere in\'egalit\'e d\'ecoulant de (\ref{systeme18}), et la seconde du lemme \ref{rang} en utilisant \'egalement le fait que $T>L$ et $T>N$. En remplacant $T$ et $N$ par leur valeur, on obtient
\begin{equation}\label{eq2003} 
h(\Delta_{\max}) \leq c_{10}L^{g-1}C_0^{(g+1)^2}(\deg_{\mathcal{L}} V\log \deg_{\mathcal{L}} V)^{g+2}(\log \log\deg_{\mathcal{L}} V)^{-2(g+2)}. 
\end{equation}

\vspace{.3cm}

 De plus, par l'in\'egalit\'e (\ref{inconnues}), et par le choix de $L$, on a 
\[ \mathfrak{b}\geq c_{11}L^{g-1}\left(L^{g+1}-c_1T(N^2)^{g-1}\deg_{\mathcal{L}} V\right)\geq c_{12}L^{(g+1)}L^{(g-1)}.	\] 
\noindent En rempla\c{c}ant $L^{g+1}$ par sa valeur, on obtient la minoration
\begin{equation}\label{eq2004} 
\mathfrak{b}\geq c_{13}C_0^{(g+1)(g+\frac{1}{2})}L^{g-1}(\deg_{\mathcal{L}} V\log \deg_{\mathcal{L}} V)^{g+1}(\log \log\deg_{\mathcal{L}} V)^{-2(g+1)}. 
\end{equation}
On reprend maintenant l'in\'egalit\'e (\ref{hauteur}) en rempla\c{c}ant les param\`etres par leurs valeurs. On obtient ainsi l'in\'egalit\'e 
\[ h(F)\leq  c_2C_0^{\frac{1}{2}(g+1)}\deg_{\mathcal{L}}(V)\log\deg_{\mathcal{L}}(V)(\log \log\deg_{\mathcal{L}} V)^{-2}.\hfill \Box \]

\vspace{.3cm}

\rem La fonction auxiliaire $F$ ainsi construite est une forme bihomog\`ene de bi-degr\'e $(L,L)$ non identiquement nulle sur $A\times A$. Elle n'est donc pas identiquement nulle sur $B$ car $N^2>L+1$.

\section{Extrapolation}
\noindent On veut montrer dans ce paragraphe que $F$ s'annule sur $i\left(\alpha_v(V)\right)$, pour $v\in \mathcal{P}_k$ appartenant \`a un ensemble convenable. Pour cela, on utilise un argument remontant \`a Dobrowolski \cite{dob} dans son c\'el\`ebre article sur la conjecture de Lehmer sur les points pour $\mathbb{G}_m$. Cet argument a \'et\'e r\'e\'ecrit et adapt\'e dans le cadre des vari\'et\'es ab\'eliennes de type C.M. dans l'article \cite{davidhindry} suivant des id\'ees de Laurent \cite{laurent}. Ce que l'on fait ici repose sur le paragraphe 6 de \cite{davidhindry}.

\begin{prop}\label{extrapolation} La fonction auxiliaire $F$ est nulle sur $i\left(\alpha_v(V)\right)$ \`a un ordre sup\'erieur \`a $T_1$ le long de $T_B$ pour toute place $v\in \mathcal{P}_k$ de norme comprise entre $\frac{1}{2}N_1$  et $N_1$.
\end{prop}
\noindent \textbf{D\'emonstration} Soit $1>\theta>\hat{\mu}^{\textnormal{ess}}_{\mathcal{L}}(V).$ Il s'agit de reprendre la proposition 6.5. de \cite{davidhindry}. On conserve donc leurs notations. Soient $v$ une place comme dans l'\'enonc\'e, $R$ un point de $V(\overline{k})$ d\'efini sur une extension $k'$ de $k$ de hauteur normalis\'ee inf\'erieure \`a $\theta$, et $w$ une place de $k'$ au dessus de $v$. Notons $\textbf{R}=(R_0,\ldots,R_n)$ un syst\`eme de coordonn\'ees projectives de $R$ dans $\mathcal{O}_w$, telles que $\mid\mid \textbf{R}\mid\mid_w=1$. Soit $\partial^{\kappa}$ un op\'erateur diff\'erentiel d'ordre $\mid\kappa\mid\leq T_1$ le long de $T_{B(\mathbb{C})}$. L'application du petit th\'eor\`eme de Fermat dans le cadre des vari\'et\'es ab\'eliennes nous donne 
\begin{equation}\label{vingt}
 \left | \partial^{\kappa}F\left(\mathbf{F}_{\alpha_v}(\mathbf{R}),\mathbf{F}^{(N)}\circ\mathbf{F}_{\alpha_v}(\mathbf{R})\right)\right |_w\leq |\pi_v|_w^{T-\mid\kappa\mid},
\end{equation}
\noindent o\`u $\mathbf{F}_{\alpha_v}$ et $\mathbf{F}^{(N)}$ sont des formes homog\`enes de $\mathcal{O}_k[\mathbf{X}]$ de degr\'e respectifs $\textnormal{N}(v)$ et $4^{m+1}$, repr\'esentant respectivement l'endomorphisme de Frobenius sur $A$ associ\'e \`a $v$, et la multiplication par $N=2^{m+1}$. Il s'agit de l'in\'egalit\'e $(20)$ p.47 de \cite{davidhindry}.

\noindent On veut maintenant sommer sur toutes les places $w$ au-dessus de $v$. Malheureusement, le choix du syst\`eme de coordonn\'ees projectives pour $R$ d\'e\-pend de $w$. On est donc obliger d'alourdir les notations pour pallier ce prob\`eme. Soient $S$, $S_N$, $S_{\alpha_v}$, $S_{N,\alpha_v}$ des coordonn\'ees projectives non nulles de $R$, $\mathbf{F}^{(N)}(\mathbf{R})$, $\mathbf{F}_{\alpha_v}(\mathbf{R})$, $\mathbf{F}^{(N)}\circ\mathbf{F}_{\alpha_v}(\mathbf{R})$ respectivement. On note de plus $S_{w,N}$, $S_{w,\alpha_v}$, $S_{w,N,\alpha_v}$ des coordonn\'ees des ces points de valeur absolue $w$-adique maximale.

\vspace{.3cm}

\noindent Soit maintenant $\partial^{\kappa}$ un op\'erateur diff\'erentiel de longueur minimale pour lequel
\[ \partial^{\kappa}F\left(\mathbf{F}_{\alpha_v}(\mathbf{R}),\mathbf{F}^{(N)}\circ\mathbf{F}_{\alpha_v}(\mathbf{R})\right)\]
\noindent est non nul. Si $|\kappa|$ est sup\'erieur \`a $T_1$, on a gagn\'e. Sinon on applique la formule de Leibniz en utilisant que $F$ est bihomog\`ene de bidegr\'e $(L,L)$. On a donc 
\begin{align*}
\partial^{\kappa}F\left(\frac{\mathbf{F}_{\alpha_v}(\mathbf{R})}{S_{\alpha_v}},\frac{\mathbf{F}^{(N)}\left(\mathbf{F}_{\alpha_v}(\mathbf{R})\right)}{S_{N,\alpha_v}}\right) & =  \frac{\partial^{\kappa}F\left(\mathbf{F}_{\alpha_v}(\mathbf{R}),\mathbf{F}^{(N)}\left(\mathbf{F}_{\alpha_v}(\mathbf{R})\right)\right)}{S_{\alpha_v}^LS_{N,\alpha_v}^L}
\end{align*}	
\noindent Or ceci est \'egal \`a 
\[ \partial^{\kappa}F\left(\frac{\mathbf{F}_{\alpha_v}(\mathbf{R})}{S_{w,\alpha_v}},\frac{\mathbf{F}^{(N)}\left(\mathbf{F}_{\alpha_v}(\mathbf{R})\right)}{S_{w,N,\alpha_v}}\right)\cdot \frac{\left(S_{w,\alpha_v}S_{w,N,\alpha_v}\right)^L}{\left(S_{\alpha_v}S_{N,\alpha_v}\right)^L}.\]
\noindent On r\'e\'ecrit alors l'in\'egalit\'e $(\ref{vingt})$ en passant au log, en sommant sur toutes les places $w$ au-dessus de $v$ et en notant $n_w$ les degr\'es locaux :
\begin{eqnarray}
\lefteqn{\sum_{w/v}n_w\log\left(\left|\partial^{\kappa}F\left(\frac{\mathbf{F}_{\alpha_v}(\mathbf{R})}{S_{\alpha_v}},\frac{\mathbf{F}^{(N)}\left(\mathbf{F}_{\alpha_v}(\mathbf{R})\right)}{S_{N,\alpha_v}}\right)\right|_w\right)}\\
& & \leq \left(T-|\kappa|\right)\sum_{w/v}n_w\log(|\pi_v|_w)+L\sum_{w/v}n_w\log\left(\frac{|S_{w,\alpha_v}S_{w,N,\alpha_v}|_w}{|S_{\alpha_v}S_{N,\alpha_v}|_w}\right)\label{eqnar}.
\end{eqnarray}
\noindent Or 
\begin{equation} \label{sys21}
\sum_{w/v}n_w\log(|\pi_v|_w)=[k':k]\log(|\pi_v|_v)\leq -[k':k]\log\left(\textnormal{N}(v)\right)
\end{equation}
\noindent De plus, on peut voir que 
\begin{align}
\sum_{w/v}n_w\log\left(\frac{|S_{w,\alpha_v}S_{w,N,\alpha_v}|_w}{|S_{\alpha_v}S_{N,\alpha_v}|_w}\right)& \leq [k':k]\left(h_{\mathcal{L}}(\alpha_v(R))+h_{\mathcal{L}}(N\alpha_v(R))\right)\\
	& \leq [k':k]\left(\textnormal{N}(v)\widehat{h}_{\mathcal{L}}(R)\!+N^2\textnormal{N}(v)\widehat{h}_{\mathcal{L}}(R)\!+c_{14}\right).\label{sys21'}
\end{align}
\noindent C'est l'in\'egalit\'e $(21)$ p. 49 de \cite{davidhindry}. En tenant compte du fait que le point $R$ est suppos\'e de hauteur (de N\'eron-Tate) inf\'erieure \`a $\theta$, et en injectant ceci dans $(\ref{eqnar})$, on obtient, en rempla\c{c}ant les param\`etres par leur valeur, l'in\'egalit\'e 
\begin{equation}\label{droite}
-\frac{1}{[k':k]}\sum_{w/v}n_w\log\left(\left|\partial^{\kappa}F\left(\frac{\mathbf{F}_{\alpha_v}(\mathbf{R})}{S_{\alpha_v}},\frac{\mathbf{F}^{(N)}\left(\mathbf{F}_{\alpha_v}(\mathbf{R})\right)}{S_{N,\alpha_v}}\right)\right|_w\right)\geq \frac{1}{2}T\log N_1.
\end{equation}
Il reste \`a majorer le membre de gauche de cette derni\`ere in\'egalit\'e. Notons $A$ ce membre de gauche. Par d\'efinition de la hauteur (absolue logarithmique) projective, on a 
\begin{align*}
A & \leq h\left(\left(\partial^{\kappa}F\left(\frac{\mathbf{F}_{\alpha_v}(\mathbf{R})}{S_{\alpha_v}},\frac{\mathbf{F}^{(N)}\left(\mathbf{F}_{\alpha_v}(\mathbf{R})\right)}{S_{N,\alpha_v}}\right)\right)^{-1}\right)\\
  &  = h\left(\partial^{\kappa}F\left(\frac{\mathbf{F}_{\alpha_v}(\mathbf{R})}{S_{\alpha_v}},\frac{\mathbf{F}^{(N)}\left(\mathbf{F}_{\alpha_v}(\mathbf{R})\right)}{S_{N,\alpha_v}}\right)\right),
\end{align*}
\noindent ceci ayant un sens grace \`a l'hypoth\`ese de non nullit\'e de $\partial^{\kappa}F(\cdots)$. Il ne reste maintenant plus qu'\`a majorer cette derni\`ere hauteur. Il s'agit d'un calcul classique (cf. par exemple \cite{davidhindry} p. 50). On obtient
\begin{equation}\label{gauche}
A\leq c_{15}\left(T_1\log(T_1+L)+LN^2\textnormal{N}(v)\theta+h(F)\right).
\end{equation}
\noindent Finalement, en mettant ensemble les in\'egalit\'es $(\ref{droite})$ et $(\ref{gauche})$, on obtient
\begin{equation}\label{extra1}
T\log N_1\leq c_{16}T_1\log(T_1+L)+c_{16}LN^2N_1\theta+c_{16}h(F).
\end{equation}
\noindent On remplace les diff\'erents param\`etres par leurs valeurs, et on obtient pour le membre de gauche de l'in\'egalit\'e, 
\[ C_0^{g+1} \deg_{\mathcal{L}}(V)\log \deg_{\mathcal{L}}(V)(\log\log \deg_{\mathcal{L}}(V))^{-2}, \]
\noindent et pour le membre de droite,
\[c_{17} C_0^g\deg_{\mathcal{L}}(V)\log \deg_{\mathcal{L}}(V)(\log\log \deg_{\mathcal{L}}(V))^{-2}.\]
D\`es que $C_0$ est assez grand, on aboutit \`a une contradiction.\hfill$\Box$

\section{Conclusion}

\noindent On commence par minorer le degr\'e de l'union des vari\'et\'es transform\'ees de $V$. 

\begin{prop}\label{somme} Soient $A$ une vari\'et\'e ab\'elienne sur $k$ de dimension $g~\geq~1$, $\mathcal{L}$ un fibr\'e en droites ample sur $A$, et $V$ une sous-$k$-vari\'et\'e stricte de $A$, irr\'eductible sur $k$.  On suppose que $V$ n'est pas une r\'eunion de sous-vari\'et\'es de torsion de $A$, et que le nombre $M$ de composantes g\'eom\'etriques de $V$ est major\'e par $c_3\deg_{\mathcal{L}}(V)^g$.  On consid\`ere enfin un ensemble d'isog\'enies $\beta_v$ admissibles deux \`a deux premi\`eres entre elles, avec $v\in\mathcal{P}_k^1=\mathcal{P}_k\cap[\![\frac{N_1}{2},N_1]\!]$. On a :
\[ \deg \left(\bigcup_{v \in \mathcal{P}_k^1}\beta_v(V)\right) \geq c_{4}\frac{\deg_{\mathcal{L}}(V)\ N_1^{g-\dim G_V}}{\log N_1}.\]
\end{prop}
\noindent \textbf{D\'emonstration} Soit $W$ une composante g\'eom\'etriquement irr\'eductible de $V$. Pour $v\in \mathcal{P}_k^1$, on a, $\beta_v$ \'etant d\'efinie sur $k$,
\[ \textnormal{card}\left(\ker(\beta_v)\cap G_{\sigma(W)}\right)=\textnormal{card}\left(\ker(\beta_v)\cap G_W\right).\]
Par ailleurs, comme $W$ n'est pas une sous-vari\'et\'e de torsion de $A$ (sinon $V$ serait r\'eunion de telles sous-vari\'et\'es), le point $(2)$ du lemme \ref{distinct} nous indique que l'\'egalit\'e $\beta_v(W)=\beta_v\left(\sigma(W)\right)$ (et $W\not=\sigma(W))$ n'est possible que pour au plus $\frac{\log M}{2}\leq c_{17}\log\deg_{\mathcal{L}}V$ \'el\'ements $v$ de $\mathcal{P}_k^1$. Notons $\mathcal{P}_k^{1\star}$ le sous-ensemble de $\mathcal{P}_k^1$ obtenu en enlevant ces \'el\'ements.  Le th\'eor\`eme de Chebotarev nous indique que 
\begin{equation}\label{chebo}
\textnormal{card}\left(\mathcal{P}_k^1\right)\geq c_{18}\frac{N_1}{\log N_1}.
\end{equation}
\noindent En rempla\c{c}ant $N_1$ par sa valeur, on constate que 
\begin{equation}\label{moinsfin}
\textnormal{card}\left(\mathcal{P}_k^{1\star}\right)\geq \frac{1}{2}\textnormal{card}\left(\mathcal{P}_k^1\right).
\end{equation}
\noindent En utilisant l'additivit\'e du degr\'e et les lemmes pr\'ec\'edents, on a 
\[ \deg_{\mathcal{L}} \left(\bigcup_{v \in \mathcal{P}_k^1}\beta_v(V)\right) \geq \deg_{\mathcal{L}} \left(\bigcup_{v \in \mathcal{P}_k^{1\star}, \sigma\in \textnormal{Gal}(\overline{k}/k)}\beta_v(\sigma(W))\right).\]
\noindent Par le lemme \ref{distinct}, ceci est sup\'erieur \`a 					
\[ \sum_{v\in\mathcal{P}_k^{1\star}}\deg_{\mathcal{L}} \left(\bigcup_{\sigma\in \textnormal{Gal}(\overline{k}/k)}\beta_v(\sigma(W))\right).\]
\noindent Enfin, le lemme \ref{degre} nous donne l'in\'egalit\'e 
\[ \deg_{\mathcal{L}} \left(\bigcup_{v \in \mathcal{P}_k^1}\beta_v(V)\right) \geq M\deg_{\mathcal{L}} W \sum_{v\in \mathcal{P}_k^{1\star}}\frac{\textnormal{q}(\beta_v)^{\dim V}}{\mid G_W\cap \textnormal{ker}(\beta_v)\mid}.\]
\noindent Le lemme 2.1. (ii) de \cite{davidhindry} nous indique que 
\[\deg_{\mathcal{L}} G_W=\left[G_W : G_W^0\right]\deg_{\mathcal{L}}(G_W^0)\leq\deg_{\mathcal{L}}(V)^g\]
\noindent En particulier on en d\'eduit que 
\[\left[G_W : G_W^0\right]\leq \deg_{\mathcal{L}}(V)^g.\]
\noindent De plus, les $\beta_v$ \'etant premiers entre eux, on a 
\begin{equation}\label{acompose}
\prod_{v\in\mathcal{P}_k^{1\star}}\mid\ker(\beta_v)\cap G_W\mid=\left| \ker\left(\prod_{v\in\mathcal{P}_k^{1\star}}\beta_v\right)\cap G_W\right|.
\end{equation}
\noindent En appliquant le lemme \ref{lemme25}, on en d\'eduit
\begin{equation}\label{acompose2}
\left| \ker\left(\prod_{v\in\mathcal{P}_k^{1\star}}\beta_v\right)\cap G_W\right|\leq \left[G_W : G_W^0\right]\left(\prod_{v\in\mathcal{P}_k^{1\star}}\textnormal{q}(\beta_v)\right)^s.
\end{equation}
\noindent En appliquant l'in\'egalit\'e arithm\'etico-g\'eom\'etrique, on obtient
\begin{equation}\label{eqfin}
 \deg \left(\bigcup_{v \in \mathcal{P}_k^1}\beta_v(V)\right)\geq c_{4}\deg_{\mathcal{L}}(V)\frac{\textnormal{card}\left(\mathcal{P}_k^{1\star}\right)N_1^{g-1}}{\left(\prod_{v\in\mathcal{P}_k^{1\star}}\mid\ker(\beta_v)\cap G_W\mid\right)^{\frac{1}{\textnormal{card}\left(\mathcal{P}_k^{1\star}\right)}}}\ \ . 
\end{equation}
\noindent En appliquant l'in\'egalit\'e (\ref{acompose2}) et la minoration du cardinal de $\mathcal{P}_k^{1\star}$, on obtient
\begin{equation}\label{eqfin2}
 \deg \left(\bigcup_{v \in \mathcal{P}_k^1}\beta_v(V)\right)\geq c_{4}\frac{\deg_{\mathcal{L}}(V)\ N_1^g}{\log N_1\left[G_W : G_W^0\right]^{\frac{\log N_1}{c_{56}N_1}}\prod_{v\in\mathcal{P}_k^{1\star}}\textnormal{q}(\beta_v)^{\frac{g}{\mathcal{P}_k^{1\star}}}}\ \ . 
\end{equation}
\noindent Enfin par d\'efinition de $\mathcal{P}_k^1$, on a la majoration $\textnormal{q}(\beta_v)\leq N_1$. En appliquant ceci et la majoration de l'indice de $G_W^0$ dans $G_W$, on a
\begin{equation}\label{eqfin3}
 \deg \left(\bigcup_{v \in \mathcal{P}_k^1}\beta_v(V)\right)\geq c_{4}\frac{\deg_{\mathcal{L}}(V)\ N_1^{g-\dim G_V}}{\log N_1}. \hspace{2cm} \Box
\end{equation}

\vspace{.3cm}

\rem C'est uniquement pour assurer l'in\'egalit\'e (\ref{moinsfin}) que l'on est conduit \`a choisir l'exposant du terme $\log \log$ dans $N_1$ tel qu'indiqu\'e, plut\^ot que l'exposant $-\frac{1}{g-\dim G_V}$ qui serait plus proche des choix de \cite{amodavid2}. Cette am\'elioration dans \cite{amodavid2} est rendue possible par la r\'esolution de deux complications techniques~: passage \`a une hypersurface secondaire explicitement construite, et raffinement galoisien.

\vspace{.3cm}

\noindent Ceci \'etant, on peut maintenant d\'emontrer le th\'eor\`eme recherch\'e. 

\vspace{.3cm}

\noindent \textbf{D\'emonstration :} on suppose par l'absurde que l'in\'egalit\'e du th\'eor\`eme \`a prouver n'est pas v\'erifi\'ee pour $C_0=c(A,\mathcal{L})^{-\frac{1}{g+3}}$ assez grand (i.e. $c(A,\mathcal{L})$ suffisamment petit). Dans cette preuve, on consid\`ere, pour all\'eger les notations, la vari\'et\'e ab\'elienne $A$, comme \'etant plong\'ee dans $\mathbb{P}_n$. Notons $\mathcal{Z}$ l'hypersurface sur $k$ de $\mathbb{P}_n$ associ\'ee \`a la forme $F\circ \varphi$ de degr\'e $(N^2+1)L$. Par choix de $N$ (\`a savoir $(N^2+1)>L$), la vari\'et\'e $\mathcal{Z}\cap A$ est une hypersurface de $A$. De plus, par la proposition \ref{extrapolation}, on sait que cette hypersurface contient les vari\'et\'es irr\'eductibles $\alpha_v(V)$ avec une multiplicit\'e sup\'erieure \`a $T_1$, pour toute place $v$ de norme comprise entre $\frac{1}{2}N_1$ et $N_1$. Donc le th\'eor\`eme de B\'ezout g\'eom\'etrique nous donne :
\[ T\deg_{\mathcal{L}}(V)+T_1\deg_{\mathcal{L}} \left(\bigcup_{\frac{N_1}{2}\leq \textnormal{N}(v)\leq N_1}\alpha_v(V)\right)\leq (\deg_{\mathcal{L}}A) L(N^2+1).\]
\noindent Cette in\'egalit\'e implique en particulier que le nombre $M$ de composantes g\'eom\'etriquement irr\'eductibles de $V$ est major\'e par une expression de la forme $c_3\deg_{\mathcal{L}}(V)^g$. On peut donc appliquer la proposition \ref{somme} avec $\beta_v=\alpha_v$. Celle-ci et l'in\'egalit\'e obtenue par le th\'eor\`eme de B\'ezout nous fournissent l'in\'egalit\'e
\begin{equation}\label{fin1}
 T_1\frac{\deg_{\mathcal{L}}(V)\ N_1^{g-\dim G_V}}{\log N_1}\leq c_{19}(A) L(N^2+1),
\end{equation}
\noindent On remplace maintenant les param\`etres par leurs valeurs pour conclure. Si $C_0$ est assez grand, l'in\'egalit\'e est contredite : si $s=g-1$, les deux membres sont du m\^eme ordre de grandeur, or, dans le membre de gauche, on a un terme constant de la forme $C_0^{2g+2}$, alors que dans le membre de droite, le terme constant est de la forme $C_0^{2g+\frac{3}{2}}$ ; sinon l'ordre de grandeur du membre de gauche est sup\'erieur \`a celui du terme de droite. (En fait, $T_1$ est construit exactement pour contredire cette in\'egalit\'e). \hfill$\Box$

\bibliographystyle{plain} \bibliography{hypersurface3}

\end{document}